\documentclass[12pt,a4paper]{article}
\usepackage{theorem}
\newtheorem{lemma}{Lemma}
\newtheorem{proposition}[lemma]{Proposition}
\newtheorem{theorem}[lemma]{Theorem}

{\theorembodyfont{\upshape}}
{\theorembodyfont{\upshape}}
{\theorembodyfont{\upshape}}
{\theorembodyfont{\upshape}\newtheorem{example}[lemma]{Example}}
{\theorembodyfont{\upshape}}
\newcommand{\N}{{\bf N}}
\newcommand{\Z}{{\bf Z}}
\newcommand{\R}{{\bf R}}
\newcommand{\C}{{\bf C}}
\newcommand{\T}{{\bf T}}
\newcommand{\rme}{{\rm e}}
\newcommand{\rmd}{{\rm d}}

\newcommand{\cB}{{\cal B}}

\newcommand{\cL}{{\cal L}}
\newcommand{\cM}{{\cal M}}
\newcommand{\cN}{{\cal N}}

\newcommand{\cW}{{\cal W}}

\newcommand{\sig}{\sigma}
\newcommand{\alp}{\alpha}

\newcommand{\lam}{\lambda}

\newcommand{\Dom}{{\rm Dom}}

\newcommand{\Spec}{{\rm Spec}}

\newcommand{\norm}{\Vert}

\newcommand{\supp}{{\rm supp}}
\newcommand{\Proof}{\underbar{Proof}{\hskip 0.1in}}

\newcommand{\lin}{{\rm lin}}
\newcommand{\clin}{\overline{{\rm lin}}}

\newcommand{\half}{\textstyle{\frac{1}{2}}}

\newcommand{\hide}[1]{}
\usepackage[final]{graphics}
\usepackage{hyperref}
\title{Triviality of the\\ Peripheral Point Spectrum}
\author{E B Davies}
\date{1 January 2005}
\begin{document}
\maketitle
\begin{abstract}
If $T_t=\rme^{Zt}$ is a positive one-parameter contraction
semigroup acting on $l^p(X)$ where $X$ is a countable set and
$1\leq p <\infty$, then the peripheral point spectrum $P$ of $Z$
cannot contain any non-zero elements. The same holds for Feller
semigroups acting on $L^p(X)$ if $X$ is locally compact.
\end{abstract}

\section{Introduction}

Let $T_t=\rme^{Zt}$ be a one-parameter contraction semigroup with
generator $Z$ acting on the separable Banach space $\cB$. We say
that $\theta\in\R$ lies in its peripheral point spectrum $P$ if
the linear subspace
\begin{eqnarray*}
\cL_\theta&=&\{f\in\Dom(Z):Zf=i\theta f\}\\
&=& \{ f\in\cB:T_tf=\rme^{i\theta t} f\mbox{ for all } t\geq 0\}.
\end{eqnarray*}
is non-zero. Under the conditions of our paper it has long been
known that $P$ is cyclic in the sense that $\alp\in P$ implies
$\alp\Z\subseteq P$. We will prove that $P$ cannot contain any
non-zero points under certain conditions. Our main results,
Theorems~\ref{9}~and~\ref{12}, are contributions to the study of
infinite-dimensional versions of the Perron-Frobenius theorems.
See \cite{EBD,nagel} for references to the literature on this
subject.

Although $P$ may be uncountable, only a countable subset of $P$
need be considered for our purposes.

\begin{lemma} \label{1}There exist sequences
$f_n\in\Dom(Z)$ and $\theta_n\in P$ such that $Zf_n=i\theta_n f_n$
and
\[
\cL:=\clin \left( \bigcup\{ \cL_\theta:\theta\in P\}\right) =
\clin\{f_n:n\in\N \}
\]
where $\clin$ stands for the closed linear span. In particular
every $f\in \bigcup\{\cL_\theta:\theta\in P\}$ is the limit of
finite linear combinations of $\{f_n:n\in\N\}$.
\end{lemma}

\Proof Since $\cB$ is separable, $\cL$ contains a countable dense
subset $\{g_n:n=1,2,...\}$. Each $g_n$ is the limit of a sequence,
each term of which is a finite linear combination of elements of
$\bigcup\{ \cL_\theta:\theta\in P\}$. The number of such elements
involved in these sums is countable.

\begin{lemma} \label{2} There
exists a sequence $t(n)\to +\infty$ such that
\begin{equation}
\lim_{n\to\infty}\norm  T_{t(n)}f-f\norm=0\label{eq1}
\end{equation}
for all $f\in\cL$.
\end{lemma}

\Proof Let $\{ f_n\} $ and $\{ \theta_n\} $ be the sequences
constructed in Lemma~\ref{1} and let $\cL_n=\lin\{ f_{r}:1\leq
r\leq n\}$. For each $n$ we will construct $t(n)\geq t(n-1)+1$
such that
\[
\norm T_{t(n)}f-f\norm \leq 2^{-n}\norm f\norm
\]
for all $f\in \cL_n$. This implies (\ref{eq1}) for all $f\in
\cL_m$ for any choice of $m$. The full statement of the lemma then
uses the fact that $T_t$ are all contractions.

Every $f\in\cL_n$ may be written in the form $f=\sum_{r=1}^n
\alp_rf_{r}$. We then have
\[
T_t f=\sum_{r=1}^n \alp_r\rme^{i\theta(r)t}f_{r}
\]
for all $t\geq 0$. Since all norms on a finite-dimensional space
are equivalent
\[
\norm (T_t-I)\vert_{\cL_n}\norm \leq c_n \max \{ |\rme^{i\theta(r)
t}-1|:1\leq r\leq n\}.
\]
The lemma follows by using the fact that the set
\[
\left\{ (\rme^{i\theta(1) t},...,\rme^{i\theta(n)
t}):t=1,2,...\right\}
\]
is a semigroup in the compact $n$-dimensional torus $\T^n$ where
$\T=\{z:|z|=1\}$. Every such semigroup has points arbitrarily
close to the group identity.

\begin{example}\label{3}
We normalize the Haar measure of $\T$ to have mass $1$, and let
$X=\T^\N$ with the usual product measure. Given $1\leq p <\infty$
define $T_t$ on $\cB=L^p(X)$ by $(T_tf)(x)=f(x_t)$, where
$x_{t,m}=\rme^{2\pi it/m}x_m$ for all $m\in\N$. Then $T_t$ is a
group of isometries on $\cB$. By taking Fourier transforms it may
be shown that $\Spec(Z)=i\R$. If we put $t(n)=n!$ then it is easy
to show that
\[
\lim_{n\to\infty} \norm T_{t(n)}f-f\norm =0
\]
for all $f\in\cB$.
\end{example}

\begin{theorem}\label{4} If
\[
\cM=\{ f\in\cB:\lim_{n\to\infty}\norm T_{t(n)}f-f\norm=0\}
\]
then $T_t$ acts as a one-parameter {\rm group} of isometries on
$\cM$.
\end{theorem}

\Proof If $f\in\cM$ and $t\geq 0$ then it follows from
\[
\norm T_{t(n)}(T_tf)-(T_tf)\norm\leq \norm T_{t(n)}f-f\norm
\]
that $T_tf\in\cM$.

If $t\geq 0$ then
\[
\norm T_{t(n)}f\norm= \norm T_{t(n)-t}T_tf\norm \leq\norm
T_tf\norm
\]
for all $t(n)\geq t$. If $f\in\cM$ then by letting $n\to\infty$ we
see that $\norm f\norm \leq\norm T_tf\norm$. Therefore $T_t$ is an
isometry when restricted to the subspace $\cM$.

It follows from the above that for any $t\geq 0$, $\cM_t=T_t(\cM)$
is a closed linear subspace of $\cM$. Moreover
$T_s(\cM_t)\subseteq\cM_t$ for all $s\geq 0$. If $f\in\cM$ then
putting $s=t(n)-t$ we deduce that
\[
f=\lim_{n\to\infty} T_{t(n)-t}(T_tf)\in\cM_t.
\]
Therefore $\cM_t=\cM$ for all $t\geq 0$. The operators $T_t$ are
defined for $t<0$ by $T_t=(T_{-t})^{-1}$.

\section{Positive Semigroups}

We now assume that $\cB=L^p_\C(X,\rmd x)$ where $1\leq p<\infty$
and $\rmd x$ is a countably additive $\sig$-finite measure on $X$.
In order to pass back and forth between the real and complex
spaces the following standard proposition is needed.

\begin{proposition}\label{5}
Let $1\leq p,q\leq\infty$ and let $A_\R :L^p_\R(X,\rmd x)\to
L^q_\R(X,\rmd x)$ be a positive linear operator with %
complex-linear extension $A_\C$. Then
\[
|A_\C ( f+i g)|\leq A_\R (|f+ig|)
\]
for all $f,g\in L^p_\R(X,\rmd x)$. Hence $\norm A_\C\norm =\norm
A_\R\norm$.
\end{proposition}

\Proof Given $\theta\in\R$ we have
\begin{eqnarray*}
|(A_\R f)\cos(\theta)+(A_\R g)\sin(\theta)|&=&
|A_\R (f\cos(\theta)+g\sin(\theta))|\\
&\leq & A_\R (|f\cos(\theta)+g\sin(\theta)|)\\
&\leq &A_\R(|f+ig|).
\end{eqnarray*}
Let $u,v,w:X\to\R$ be functions in the classes of $A_\R f,A_\R
g,A_\R (|f+ig|)$. Then we have shown that
\[
|u(x)\cos(\theta)+v(x)\sin(\theta)|\leq w(x)
\]
for all $x$ not in some null set $N(\theta)$. If $\{
\theta_n\}_{n=1}^\infty$ is a countable dense subset of
$[-\pi,\pi]$ then
\[
|u(x)+iv(x)|=\sup_{1\leq n
<\infty}|u(x)\cos(\theta_n)+v(x)\sin(\theta_n)|\leq w(x)
\]
for all $x$ not in the null set $\bigcup_{n=1}^\infty
N(\theta_n)$. This implies the first statement of the theorem,
from which the second follows immediately.

We now assume that $T_t$ is a contraction semigroup on $\cB$ which
commutes with complex conjugation and which is positive when
restricted to $L^p_\R(X,\rmd x)$. This latter space is a lattice
with respect to the operations
\begin{eqnarray}
f\vee g&=&\max\{f,g\}=\half(f+g)+\half|f-g|,\label{eq2}\\
f\wedge g&=&\min\{f,g\}=\half(f+g)-\half|f-g|.\label{eq3}
\end{eqnarray}

The strict monotonicity of the norm in $L^p_\R(X,\rmd x)$ for
$1\leq p<\infty$ is of critical importance for our next theorem
and hence for the rest of the paper.

\begin{theorem}\label{6}
The subspace $\cM$ is closed under complex conjugation and its
real part $\cN$ is a closed linear sublattice of $L^p_\R(X,\rmd
x)$. Moreover the operators $T_t$ are lattice isomorphisms when
restricted to $\cN$ for all $t\in\R$.
\end{theorem}

\Proof We use the fact, proved in Proposition~\ref{5}, that
$|T_t(f)|\leq T_t(|f|)$ for all $f\in L^p_\C(X,\rmd x)$ and $t\geq
0$. If $f\in\cM$ and $t\geq 0$ then
\[
\norm \,|f|\,\norm\geq \norm T_t(|f|)\norm\geq
\norm\,|T_t(f)|\,\norm=\norm T_tf\norm=\norm
f\norm=\norm\,|f|\,\norm.
\]
It follows that $\norm T_t(|f|)\norm= \norm\,|T_t(f)|\,\norm$, and
combining this with $|T_t(f)|\leq T_t(|f|)$ and $1\leq p <\infty$,
we deduce that $|T_t(f)|= T_t(|f|)$. This finally implies that
\[
\norm T_{t(n)}(|f|)-|f|\,\norm =\norm
\,|T_{t(n)}(f)|-|f|\,\norm\leq \norm T_{t(n)}(f)-f\norm\to 0
\]
as $n\to\infty$. Therefore $|f|\in\cM$.

This establishes that $\cN$ is a closed linear sublattice of
$L^p_\R(X,\rmd x)$. The identity $T_t(|f|)=|T_t(f)|$ for all
$f\in\cM$ and $t\geq 0$ implies the same equality for all
$t\in\R$. Hence $T_t$ restricted to $\cN$ is a lattice isomorphism
for all $t\in\R$ by (\ref{eq2}) and (\ref{eq3}).

\begin{lemma}\label{7}
There exists a Borel set $E$ in $X$ and $e\in\cN_+$ such that
$\supp(e)=E$ and $\supp(g)\subseteq E$ for all $g\in\cN$. In
particular $\supp(f_\theta)\subseteq E$ for all $\theta \in P$.
The set $E$ is invariant in the sense that if $g\in\cB$ and
$\supp(g)\subseteq E$ then $\supp(T_t(g))\subseteq E$ for all
$t\geq 0$.
\end{lemma}

\Proof Let $e_n$ be a countable dense subset of $\cN$ and let $E$
be the support of
\[
e:=\sum_{n=1}^\infty 2^{-n}|e_n|/\norm e_n\norm.
\]
If $g\in\cN$ then by approximating $g$ by terms of the sequence
$e_n$ we see that $\supp(g)\subseteq E$. In particular
$\supp(T_te)\subseteq E$ for all $t\geq 0$. This implies that $E$
is an invariant set.

Our next proposition, which is an $l^p$ version of the
Stone-Weierstrass theorem, depends on making some further
definitions. Let $\cW$ be any closed linear sublattice of
$l^p(X,w)$ where $1\leq p <\infty$ and $w(x)$ is a positive weight
for every $x$ in the countable set $X$ and $\sum_{x\in X}
w(x)<\infty$. Suppose that $1\in\cW$. We define an equivalence
relation on $X$ by putting $x\sim y$ if $f(x)=f(y)$ for all
$f\in\cW$ and let $\tilde{X}$ denote the set of all equivalence
classes.

\begin{proposition} \label{8}
If $E\in\tilde{X}$ then its characteristic function $\chi_E$ lies
in $\cW$. Moreover
\[
\cW=\clin\{ \chi_E:E\in\tilde{X} \}.
\]
\end{proposition}

\Proof If $f\in\cW$ then $f$ is constant on any class $u$, and we
write $f(u)$ to denote this constant value.

Let $v$ be some chosen class. If $u$ is a distinct class then
there exists $g_{u}\in\cW$ such that $g_{u}(u)> g_{u}(v)$. On
putting $h_{u}=(g_{u}-g_{u}(v)1)\vee 0$ we see that
$h_{u}\in\cW_+$, $h_{u}(u)>0$ and $h_{u}(v)=0$. If $n(\cdot)$ is
some enumeration of $\tilde{X}$ then
\[
k=\sum_{\{u:u\not= v\}}2^{-n(u)} h_{u}/\norm h_{u}\norm_p
\]
lies in $\cW_+$ and satisfies $k(v)=0$ and $k(u)>0$ for all
$u\not= v$. If $p_{n}=1\wedge (nk)$ then $p_{n}$ increases
monotonically with $n$ and converges in $l^p$ norm to $1-\chi_v$
by the dominated convergence theorem. Hence $\chi_v\in \cW$.

The final statement of the proposition follows from the fact that
if $f\in\cW_+$ then
\[
f=\sum_{u\in\tilde{X}} f(u)\chi_u
\]
where this sum converges monotonically and in $l^p$ norm by virtue
of the dominated convergence theorem.

\begin{theorem}\label{9}
If $T_t$ is a positive one-parameter contraction semigroup acting
on $l^p(X)$ where $X$ is a countable set, then $T_t\vert_\cN=I$.
In particular the peripheral point spectrum $P$ cannot contain any
non-zero elements.
\end{theorem}

\Proof We henceforth assume that the set $E$ of Lemma~\ref{7}
equals $X$; equivalently we restrict attention to the invariant
subspace $l^p(E)$. We define the positive weight $w$ on $X$ by
$w(x)=e(x)^p$ where $e$ is the function defined in Lemma~\ref{7}.
We then define the isometry $U:l^p(X,w)\to l^p(X)$ by $Ug=eg$ and
transfer $T_t$ to $l^p(X,w)$ by putting $\tilde{T}_t =U^{-1}T_tU$.
On suppressing the tilde, this enables us to assume that $1\in
\cN$ and that $X$ has finite measure.

We next apply Proposition~\ref{8}. We say that $h\in\cN_+$ lies on
an extreme ray of $\cN_+$ if $0\leq g\leq h$ and $g\in\cN_+$
implies that $g=\lam h$ for some $\lam$ such that $0\leq \lam\leq
1$.  It is immediate that $h$ lies on an extreme ray if and only
if $h=\lam\chi_E$ where $\lam\geq 0$ and $E\in\tilde{X}$. Since
$T_t$ are lattice automorphisms they permute the extreme rays of
$\cN_+$ and hence induce permutations $\pi(t)$ of $\tilde{X}$.

Given an equivalence class $E$, we see that $T_t(\chi_E)$ depends
continuously on $t$, but it is also a multiple of
$\chi_{\pi(t)E}$. Therefore $T_t(\chi_E)=\chi_E$ for all $t\in\R$
and $T_t=I$ on $\cN$. This implies that $T_t=I$ on $\cM$ and hence
that $\theta=0$ for all $\theta\in P$.

\section{Feller Semigroups}

In this section we suppose that $X$ is a locally compact Hausdorff
space with a countable basis to its topology, and that $\rmd x$ is
a regular Borel measure on $X$ with support equal to $X$. We
suppose that $T_t=\rme^{Zt}$ is a positive one-parameter semigroup
on $\cB=L^p(X,\rmd x)$, where $1\leq p<\infty$. We finally suppose
that $T_t$ has the Feller property: that is $T_tf\in C(X)$ for all
$f\in\cB$ and $t>0$, where $C(X)$ is the space of all continuous
functions on $X$. Although elements of $L^p(X,\rmd x)$ are only
defined up to null sets, our support hypothesis implies that if an
element of $L^p(X,\rmd x)$ can be represented by a continuous
function, then that continuous representative is unique. If $X$ is
discrete and $\rmd x$ is the counting measure then the Feller
property is automatic, so the theorems of this section contain the
earlier results as special cases.

It follows from our assumptions that any eigenfunction of $Z$ lies
in $C(X)$. Theorem~\ref{4} states that $T_t(\cM)=\cM$ for all
$t>0$, so $\cM$ is both a closed subspace of $L^p(X,\rmd x)$ and
also a subspace of $C(X)$. The closed graph theorem implies that
within $\cM$ convergence in $L^p$ norm implies locally uniform
convergence. If $f\in C(X)$ we put
\[
\supp(f)=\{ x: f(x)\not= 0\},
\]
this being an open set.

\begin{lemma} \label{10}The real part $\cN$ of $\cM$ is a linear sublattice
of $C(X)$. There exists an invariant open set $U\subseteq X$ and
$e\in\cN_+$ such that $\supp(e)=U$ and $\supp(g)\subseteq U$ for
all $g\in\cN$.
\end{lemma}

\Proof This is a result of combining the above observations with
Theorem~\ref{6} and Lemma~\ref{7}.

The set $U$ is empty if and only if $\cM=0$, in which case the
peripheral point spectrum is empty; we assume that this is not the
case. From this point onwards we restrict attention to the action
of $T_t$ on the invariant subspace $L^p(U,\rmd x)$, which contains
all of the eigenfunctions associated with the peripheral point
spectrum $P$. This is most conveniently done by putting $X=U$.

We define the isometry $V:L^p(X,e(x)^p\rmd x)\to L^p(X,\rmd x))$
by $Vg=eg$ where $e$ is the positive continuous function of
Lemma~\ref{10}, and transfer $T_t$ to $L^p(X,e(x)^p\rmd x)$ by
putting $\tilde{T}_t =V^{-1}T_tV$. On suppressing the tilde, this
enables us to assume that $1\in \cN$ and that $X$ has finite
measure.

As before we define an equivalence relation on $X$ by putting
$x\sim y$ if $f(x)=f(y)$ for all $f\in\cN$.

\begin{lemma}\label{11} Every equivalence class $E$ is both open
and closed and has positive measure. Moreover $\chi_E\in\cN$.
\end{lemma}

\Proof We modify the proof of Proposition~\ref{8}. If $a\in E$
then
\[
E=\bigcap_{f\in \cN} \{x:f(x)=f(a)\}
\]
so $E$ is a closed set. If $x\notin E$ then there exists $g_x\in
\cN$ such that $g_x(x)>g_x(a)$. Putting $h_x=(g_x-g_x(a)1)\vee
0\in \cN_+$ and $U_x=\{y:h_x(y)>0\}$ we see that
\[
\bigcup_{x\notin E} U_x=X\backslash E.
\]
The Lindel\"of property states that there exists a sequence
$x(n)\notin E$ such that
\[
\bigcup_{n=1}^\infty U_{x(n)}=X\backslash E.
\]
It follows that if
\begin{equation}
k=\left\{\sum_{n=1}^\infty 2^{-n} h_{x(n)}/\norm
h_{x(n)}\norm_p\in\cN_+\right\} \wedge 1\label{eq4}
\end{equation}
then $k(x)=0$ for all $x\in E$ and $0 < k(x)\leq 1$ for all
$x\notin E$. Finally $p_n=1\wedge(nk)$ increases monotonically as
$n$ increases and
\begin{equation}
\lim_{n=1}p_n=1-\chi_E.\label{eq5}
\end{equation}
Therefore $\chi_E\in \cN$.

The limits in (\ref{eq4}) and (\ref{eq5}) were taken in $L^p$
norm. However, they involve monotone sequences and converge
pointwise (everywhere, not almost everywhere) to the stated
limits. In particular $k=\lim_{m\to\infty}k_m$ where
\[
k_m=\left\{\sum_{n=1}^m 2^{-n} h_{x(n)}/\norm
h_{x(n)}\norm_p\right\} \wedge 1\in\cN_+
\]
increase monotonically and are bounded above by $1$. Since all of
the functions lie in $\cN\subseteq C(X)$, $L^p$ norm convergence
implies locally uniform convergence. We deduce that $\chi_E$ is
continuous and hence that $E$ is open. Since $E$ is non-empty its
measure must be positive.

\begin{theorem}\label{12} If $T_t$ is a one-parameter positive
semigroup with the Feller property acting on $L^p(X,\rmd x)$ for
some $1\leq p<\infty$ then $T_t\vert_\cN=1$. Hence the peripheral
point spectrum $P$ cannot contain any non-zero points.
\end{theorem}

\Proof Without loss of generality we may restrict to the open set
$U$ and apply the isometry $V$ as described before Lemma~\ref{11}.
The remainder of the proof follows the argument of
Theorem~\ref{9}.

\section{Irreducibility}

Let $T_t$ be a positive one-parameter contraction semigroup acting
on $\cB=L^p_\R(X,\rmd x)$ where $1\leq p<\infty$ and $\rmd x$ is a
countably additive $\sig$-finite measure on $X$.

\begin{theorem}\label{13} Suppose that $f\in \cB_+$ has support
$S$ and $T_tf=f$ for all $t>0$. If $E\subseteq X$ is an invariant
Borel set for the semigroup $T_t$ then so is $S\backslash E$.
\end{theorem}

\Proof Since $S$ is an invariant set we may restrict attention to
the action of $T_t$ in the subspace $L^p(S,\rmd x)$. Putting
$g_t=T_t(\chi_E f)$ we observe that $0\leq g_t\leq f$ and
$\supp(g_t)\subseteq E$. Hence $0\leq g_t\leq g_0$ for all $t\geq
0$. Putting $h_t=T_t(\chi_{S\backslash E}f)$ we deduce that
\[
h_t=T_t(f-g_0)\geq f-g_0=h_0.
\]
Sine $T_t$ are contractions and the norm is strictly monotone we
deduce that $h_t=h_0$. This implies that $\supp(h_0)$, i.e.
$S\backslash E$, is invariant.

\begin{theorem}\label{14} Let $X$ be a locally compact Hausdorff
space with a countable basis to its topology, and let $\rmd x$ be
a regular Borel measure on $X$ with support equal to $X$. Let
$T_t$ be a positive one-parameter semigroup with the Feller
property acting on $\cB=L^p(X,\rmd x)$, where $1\leq p<\infty$.
Let $f$ be a positive continuous function in $\cB$ satisfying
$T_tf=f$ for all $t>0$. If $X$ is connected then $T_t$ is
irreducible.
\end{theorem}

\Proof If $E$ is an invariant Borel set in $X$ and $g=\chi_Ef$
then the proof of Theorem~\ref{13} establishes that $T_tg=g$ for
all $t>0$. If we let $g$ be the unique continuous function in the
class of the `same' element of $L^p$ then $g^2=fg$ almost
everywhere implies $g^2=fg$ everywhere. The continuous function
$h=g/f$ satisfies $h^2=h$ so the support $F$ of $h$ is both open
and closed. Since $X$ is connected either $F=X$ or $F=\emptyset$.
We finally observe that $E=F$ up to alteration on a null set.

Department of Mathematics\\
King's College\\
Strand\\
London\\
WC2R 2LS\\
UK

E.Brian.Davies@kcl.ac.uk

\end{document}